%% file: puget.tex
\documentclass[12pt]{article}

\input{macros_angl.tex}

\title{On strictly convex subsets in negatively curved manifolds} 
\author{Jouni Parkkonen \and Fr\'ed\'eric Paulin}
\date{\today}

\begin{document}
\maketitle

\begin{abstract}
  In a complete simply connected Riemannian manifold $X$ of pinched
  negative curvature, we give a sharp criterion for a subset $C$ to be
  the $\epsilon$-neighbourhood of some convex subset of $X$, in terms
  of the extrinsic curvatures of the boundary of $C$. \footnote{ 
    convexity, negative curvature, second fundamental form, Riemannian
    convolution smoothing.  }
\end{abstract}

\section{Introduction}
\label{sec:intro}

Let $X$ be a complete simply connected smooth Riemannian manifold of
dimension $m\geq 2$ with pinched negative sectional curvature
$-b^2\leq K\leq -a^2<0$, with $a,b>0$. Fix $\epsilon>0$.  A subset $C$
of $X$ will be called {\it $\epsilon$-strictly convex} if there exists
a nonempty convex subset $C'$ in $X$ such that $C$ is the closed
$\epsilon$-neighbourhood of $C'$.  One often encounters
$\epsilon$-strictly convex subsets in the litterature, for instance
when considering tubular neighbourhoods of geodesic lines or totally
geodesic subspaces (see for instance \cite{Gra}), or the
$\epsilon$-neighbourhood of the convex hull of the limit set of a
nonelementary discrete subgroup of isometries of $X$ (see for instance
\cite{MT}). In \cite{PP}, we studied penetration properties of
geodesic lines in $\epsilon$-strictly convex subsets of $X$ (called
$\epsilon$-convex subsets in \cite{PP}); given an appropriate family
of them, we constructed geodesic lines having an exactly prescribed
(big enough) penetration in exactly one of them, and otherwise
avoiding (or not entering too much in) them.

In this paper, we give a criterion for a subset of $X$ to be
$\epsilon$-strictly convex, in terms of the extrinsic curvature
properties of its boundary. Note that an $\epsilon$-strictly convex
subset of $X$ is a closed strictly convex subset of $X$ with nonempty
interior, whose boundary is  ${\rm C}^{1,1}$-smooth (see for instance
\cite[Theo.~4.8 (9)]{Fed} or \cite[page 272]{Wal}).

For every strictly convex ${\rm C}^{1,1}$-smooth hypersurface $S$ in
$X$, let ${\rm II}_S:TS\oplus TS\ra\RR$ be the (almost everywhere
defined, scalar) second fundamental form of $S$ associated to the
inward normal unit vector field $\vec{n}$ along $S$ (see Section
\ref{sect:proofmain} for a definition). For every point $x$ of $S$ at
which ${\rm II}_S$ is defined, let $\overline{\rm II}_S(x)$
(resp.~$\underline{\rm II}_S(x)$) be the least upper (resp.~greatest
lower) bound of ${\rm II}_S(v,v)$ for every unit tangent vector $v$ to
$S$ at $x$.

\btheo \label{theo:maintheo}
Let $C$ be a closed strictly convex subset of $X$ with nonempty
interior, whose boundary $S=\partial C$ is ${\rm C}^{1,1}$-smooth.
\begin{itemize}
\item[$\bullet$] 
  If $C$ is $\epsilon$-strictly convex, then 
$a\tanh (a\epsilon)\le\underline{\rm II}_S\le\overline{\rm II}_S\leq
  b\coth(b\epsilon)$ 
  almost everywhere.
\item[$\bullet$] 
  If $b\tanh (b\epsilon)\le\underline{\rm
    II}_S\le\overline{\rm II}_S\leq a\coth(a\epsilon)$ 
almost everywhere, then $C$ is
  $\epsilon$-strictly convex.
\end{itemize}
\etheo

The bounds that appear in the statement of the theorem are
geometrically natural: $a\tanh(a\epsilon)$ is the value on unit
tangent vectors of the second fundamental form of the
$\epsilon$-neighbourhood of a totally geodesic hyperplane in the real
hyperbolic $m$-space of constant curvature $-a^2$, and
$b\coth(b\epsilon)$ is the corresponding number for a sphere of radius
$\epsilon$ in the real hyperbolic $m$-space with constant curvature
$-b^2$.  Furthermore, in the complex hyperbolic plane whose sectional
curvatures are normalized to be between $-4$ and $-1$, a geodesic line
is contained in a copy of the real hyperbolic plane of curvature $-1$
and it is orthogonal to a copy of the real hyperbolic plane of
curvature $-4$. Thus, the $\epsilon$-neighbourhood of a geodesic line
has principal curvatures corresponding to both extremes
$\tanh\epsilon$ and $2\coth\epsilon$ in the first part of the
statement of the theorem.

One should expect not only upper bounds on the extrinsic curvatures
but also lower bounds.
Indeed, if for instance $C$ is a half-space in the real hyperbolic
$m$-space, then $C$ is closed, convex, with nonempty interior and
smooth boundary, but is not the closed $\epsilon$-neighbourhood of a
convex subset (the natural subset of $X$ of which it is the closed
$\epsilon$-neighbourhood is nonconvex, and even strictly concave).

A horoball in $X$ is $\epsilon$-strictly convex for every
$\epsilon>0$, hence we recover from the first assertion of Theorem
\ref{theo:maintheo} that the extrinsic curvatures of a horosphere in
$X$ belong to $[a,b]$.

Theorem \ref{theo:maintheo} is sharp, given the pinching on the
curvature, since it immediately implies the following result.

\bcoro Let $X$ be a complete simply connected smooth Riemannian
manifold of dimension $m\geq 2$ with constant sectional curvature
$-a^2<0$, let $C$ be a closed strictly convex subset of $X$ with
nonempty interior, whose boundary $S=\partial C$ is ${\rm
  C}^{1,1}$-smooth, and let $\epsilon>0$. Then $C$ is
$\epsilon$-strictly convex if and only if $a\tanh (a\epsilon)\leq
\underline{\rm II}_S\leq \overline{\rm II}_S\leq a\coth(a\epsilon)$
almost everywhere.  \cqfd 
\ecoro

Our notion of $\epsilon$-strictly convexity is related to, but
different from, the notion of $\lambda$-convexity studied for instance
in \cite{GR,BGR,BM}. Indeed, for every $\lambda\in[0,1]$, there exist
$\lambda$-convex (in the sense of these references) subsets of the
real hyperbolic plane (with constant sectional curvature $-1$) that
are not $\epsilon$-strictly convex for any $\epsilon>0$ (for instance
the intersection of all $\lambda$-convex subsets containing two
distinct points). Furthermore, the lower bound in the first assertion
of Theorem \ref{theo:maintheo} says that any $\epsilon$-strictly
convex subset with $\operatorname{C}^2$-smooth boundary in $X$ is
$a\tanh(a\epsilon)$-convex in the sense of \cite[Def.~2.2]{BGR}.

Alexander and Bishop \cite{AB} (see also
\cite{Lytchak}), have introduced a natural notion of an ``extrinsic
curvature bounded from above'' for subspaces of ${\rm
  CAT}(\kappa)$-spaces, extending the notion of having a bounded
(absolute value of the) second fundamental form for submanifolds of
Riemannian manifolds.  Thus, this concept of \cite{AB} is related to
our notion of $\epsilon$-strictly convex subsets (see in particular
Proposition 6.1 in \cite{AB}).

\medskip Comparison techniques in Riemannian geometry (as in
\cite{CE,Pet,Esch}) are at the heart of the proof of Theorem
\ref{theo:maintheo}. We start by developing the (quite classical)
vector space version of them in Section \ref{sec:compaRiccati}, in a
symmetric way to treat upper and lower bounds. The motivations for
such a vector space study will be given at the beginning of Section
\ref{sect:proofmain}, along with the definitions of the Riemannian
geometry tools that we will need. We then prove Theorem
\ref{theo:maintheo}, using the Riemannian convolution smoothing
process of Greene and Wu \cite{GW1} to deal with regularity issues.

\medskip\noindent{\small {\it Acknowledgement. } We thank P.~Pansu for
  useful conversations.}

\section{Comparison results for the matrix Riccati equation}
\label{sec:compaRiccati}

The arguments of the following proposition are standard (compare for
instance with \cite[\S 2]{Esch}, \cite[\S 6.5]{Pet}), 
but do not seem to appear in this precise form in the litterature.
Recall that an endomorphism $f$ of a (finite dimensional) Euclidean
space $E$ is said to be {\it nonnegative}, and we then write $f\geq
0$, if $\langle f(v)\,|\,v\rangle\geq 0$ for every $v\in E$, and {\em
  positive} if $\langle f(v)\,|\,v\rangle>0$ for all $v\ne0$. We
denote by $\operatorname{Sym}(E)$ the (finite dimensional) vector
space of symmetric endomorphisms of $E$.

\bprop \label{prop:matrixRiccati} Let $\epsilon,a,b>0$, let
$R\colon]-\epsilon,+\epsilon[\,\ra\operatorname{Sym}(E)$ be a smooth
map such that
\begin{equation}\label{eq:inegmatrixriccatti}
\forall\;t\in\;]-\epsilon,+\epsilon[\,,\;\;\;
a^2{\rm Id}\leq R(t) \leq b^2{\rm Id}\;,
\end{equation} 
and let $t\mapsto A(t)$ be a smooth map from a maximal neighbourhood
of $0$ in $\RR$ to $\operatorname{Sym}(E)$ such that $-\dot{A}(t) +
A(t)^2+R(t)=0$ for all $t$ in the domain of $A$.
Let $\lambda_-(t)$ and $\lambda_+(t)$ denote respectively the
smallest and biggest eigenvalues of $A(t)$. 

(i) If 
$$
b\tanh (b\epsilon)\le\lambda_-(0)\le\lambda_+(0)\le a\coth
(a\epsilon),
$$ 
then the map $t\mapsto A(t)$ is defined at least on $[0,\epsilon\,[
\,$, and $A(t)$ is positive for every $t$ in $[0,\epsilon\,[$.

(ii) Conversely, if $A(t)$ is defined and nonnegative on
$[0,\epsilon\,[$, then
$$
a\tanh (a\epsilon)\le\lambda_-(0)\le\lambda_+(0)\le b\coth (b\epsilon)\;.
$$
\eprop

\dem Let us first prove that the functions $\lambda_\pm$ are locally
Lipschitz, hence are almost everywhere differentiable, and that they
satisfy almost everywhere the inequalities
\begin{equation}\label{eq:riccatibounds}
a^2\leq -\dot\lambda_\pm(t)+\lambda_\pm^2(t)\leq b^2\;.
\end{equation}

To prove the locally Lipschitz regularity of the path of biggest
(resp.~smallest) eigenvalues $t\mapsto \sigma(t)$ (resp.~$t\mapsto
\iota(t)$) of a smooth path of symmetric endomorphisms $t\mapsto
\A(t)$ of a Euclidean space, up to adding a constant big multiple of
the identity to $\A(t)$, we may assume that $\A(t)$ is positive. Then,
using the standard operator norm, $\sigma(t)= \|\A(t)\|$ and
$\iota(t)=\|\A(t)^{-1}\|^{-1}$. Hence $\sigma$ and $\iota$ are indeed
locally Lipschitz.

To prove the inequalities \eqref{eq:riccatibounds}, let $t$ be a time
at which $\lambda'_\pm$ is defined, and let $v_t$ be a unit
eigenvector of $A(t)$ corresponding to $\lambda_+(t)$.  In particular,
$\lambda_+(t)=\langle A(t)v_t\,|\,v_t\rangle$ and
$\lambda_+^2(t)=\langle A^2(t)v_t\,|\, v_t\rangle$. For every $\eta>0$
small enough, since by definition,  $\lambda_+(t-\eta)=\max_{\|w\|=1}\langle
A(t-\eta)w\,|\,w\rangle$, we have
\begin{equation}\label{eq:upperderiv}
\langle\big(A(t-\eta)-A(t)\big)v_t\,|\,v_t\rangle\leq  
\lambda_+(t-\eta)-\lambda_+(t)\;.
\end{equation}
Dividing by $-\eta<0$ and taking the limit as $\eta$ goes to $0$ gives
$\langle A'(t)v_t\,|\,v_t\rangle\geq \lambda_+'(t)$.  It follows that
$$
-\dot\lambda_+(t)+\lambda_+^2(t)\geq -\langle A'(t)v_t\,|\,v_t\rangle+ 
\langle A^2(t)v_t\,|\, v_t\rangle\geq a^2\;.
$$
The inequality $-\dot\lambda_+(t)+\lambda_+^2(t)-b^2\leq 0$ is proved
similarly, by replacing $\eta$ by $-\eta$ in the formula
\eqref{eq:upperderiv}, and dividing by $\eta>0$. The inequalities
$a^2\leq -\dot\lambda_-(t)+\lambda_-^2(t) \leq b^2$ are also proved
similarly.

\medskip 
The proposition \ref{prop:matrixRiccati} will follow from the
next two lemmae.

\blemm \label{lem4}
Let $\epsilon>0$. Let $s:t\mapsto s(t)$ be a real locally
Lipschitz map defined on an open interval $I$ of $\RR$ containing $0$,
such that $-s'(t)+s^2(t)-a^2\geq 0$ almost everywhere and $s(0)\leq
a\coth (a\epsilon)$. Then $s(t)\leq a\coth a(\epsilon-t)$ for every
$t\in I\cap[0,\epsilon[\,$. Conversely, if $s$ is a real locally
Lipschitz map defined on $[0,\epsilon[$ such that
$-s'(t)+s^2(t)-b^2\leq 0$ almost everywhere, then $s(0)\leq b\coth
(b\epsilon)$.  
\elemm

\dem For every $c>0$, the maximal solution $s_{c,\epsilon}$ of the
scalar differential equation $-x'+x^2-c^2=0$, with value $c\coth
(c\epsilon)$ at $t=0$, is $s_{c,\epsilon}:t\mapsto c\coth
c(\epsilon-t)$, defined on $]-\infty, \epsilon[$, which satisfies
$\lim_{t\ra\epsilon^-}s_{c,\epsilon}(t)=+\infty$.  If
$$
\varphi_{c,\epsilon}(t)= \big(s_{c,\epsilon}(t)-s(t)\big) e^{-\int_0^t
  (s_{c,\epsilon}(u)+s(u))\,du},
$$ 
which is defined on the intersection of the intervals of definition of
$s_{c,\epsilon}$ and $s$, then
$$
\varphi'_{c,\epsilon}(t) =
\big(s'_{c,\epsilon}(t)- s_{c,\epsilon}^2(t) -s'(t)+s^2(t)\big)
e^{-\int_0^t(s_{c,\epsilon}(u)+s(u))\,du}
$$ 
almost everywhere.

If $-s'(t)+s^2(t)-a^2\geq 0$ almost everywhere, then
$\varphi'_{a,\epsilon} \geq 0$ almost everywhere. Hence
$\varphi_{a,\epsilon}$ is nondecreasing. If furthermore $s(0)\leq
a\coth (a\epsilon)$, then $\varphi_{a,\epsilon}(0) =s_{a,\epsilon}(0)
-s(0)\geq 0$, so that the map $\varphi_{a,\epsilon}$ is nonnegative at
nonnegative times. Hence, for every $t\geq 0$, we have $s(t)\leq
s_{a,\epsilon}(t)$ whenever both functions are defined. The
first claim follows.

To prove the second one, if by absurd $s(0)> b\coth (b\epsilon)$,
since the map $\epsilon\mapsto b\coth (b\epsilon)$ is continuous,
there exists $\epsilon'<\epsilon$ such that $s(0)\geq b\coth
(b\epsilon')$.  If $-s'(t)+s^2(t)\leq b^2$ almost everywhere, then the
map $\varphi_{b,\epsilon'}$ is nonincreasing. Since
$\varphi_{b,\epsilon'}(0)\leq 0$, $\varphi_{b,\epsilon'}$ is
nonpositive at nonnegative times. Therefore $s(t)\geq
s_{b,\epsilon'}(t)$ for $t\in[0,\epsilon'[\,$, so that $s(t)$ tends to
$+\infty$ as $t\ra {\epsilon'}^-$. But this contradicts the finiteness
of $s(\epsilon')$, since $s$ in continuous at $\epsilon'$. 
\cqfd

\blemm \label{lem5}
Let $\epsilon>0$. Let $i:t\mapsto i(t)$ be a real locally
Lipschitz map defined on $[0,\epsilon[$. If $-i'(t)+i^2(t)-b^2\leq 0$
almost everywhere and $i(0)\geq b\tanh (b\epsilon)$, then $i$ is
positive on $[0,\epsilon[\,$.  Conversely, if $i$ is nonnegative on
$[0,\epsilon[$ and $-i'(t)+i^2(t)-a^2\geq 0$ almost everywhere, then
$i(0)\geq a\tanh (a\epsilon)$.  
\elemm

\dem 
The proof is similar to that of the previous lemma. For every $c>0$,
the maximal solution $i_{c,\epsilon}$ of the scalar differential
equation $-x'+x^2-c^2=0$, with value $c\tanh (c\epsilon)$ at $t=0$, is
$i_{c,\epsilon}:t\mapsto c\tanh c(\epsilon-t)$, defined on $\RR$.  If
$$
\psi_{c,\epsilon}(t)=\big(i(t)-i_{c,\epsilon}(t)\big)
e^{-\int_0^t(i(u)+i_{c,\epsilon}(u))\,du},
$$ 
then 
$$
\psi'_{c,\epsilon}(t)= \big(i'(t)-i^2(t)
-i'_{c,\epsilon}(t)+i_{c,\epsilon}^2(t)\big)
e^{-\int_0^t(i(u)+i_{c,\epsilon}(u))\,du}
$$ 
almost everywhere.

If $-i'(t)+i^2(t)-b^2\leq 0$ almost everywhere, then
$\psi'_{b,\epsilon}(t) \geq 0$ almost everywhere. Hence
$\psi_{b,\epsilon}$ is nondecreasing. If furthermore $i(0)\geq b\tanh
(b\epsilon)$, then $\psi_{b,\epsilon}(0)=i(0)-i_{b,\epsilon}(0)\geq
0$, so that the map $\psi_{b,\epsilon}$ is nonnegative at nonnegative
times. Hence, if $t\geq 0$, we have $i(t)\geq i_{b,\epsilon}(t)$
whenever defined. Since $i_{b,\epsilon}(t)> 0$ if and only if
$t<\epsilon$, the first assertion follows.

To prove the second one, assume that $-i'(t)+i^2(t)-a^2\geq 0$ almost
everywhere. Then $\psi_{a,\epsilon}$ is nonincreasing. If furthermore
$i(t)\geq 0$ for $t\in [0,\epsilon[\,$, then $\liminf_{t\ra\epsilon^-}
\psi_{a,\epsilon}(t)\geq 0$ since $i_{a,\epsilon}(\epsilon)=0$, so
that $\psi_{a,\epsilon}(0)\geq 0$, which implies that $i(0)\geq a\tanh
(a\epsilon)$.  
\cqfd

\medskip Let us now prove Assertion (i) of Proposition
\ref{prop:matrixRiccati}. The first claim of Lemma \ref{lem4} applied
(using Equation \eqref{eq:riccatibounds}) to the biggest eigenvalue
$s(t)=\lambda_+(t)$ of $A(t)$ shows that the matrix $A(t)$ remains
bounded as long as it is defined and $t\leq\epsilon$. By the
assumption of maximality on the domain of definition of $A$, and by
the non explosion theorem at a finite bound for ordinary differential
equations, this proves that $A(t)$ is defined for every $t$ in the
interval $[0,\epsilon[\,$. The fact that $A(t)$ is positive there
follows from the first claim of Lemma \ref{lem5} applied (using
Equation \eqref{eq:riccatibounds}) to $i=\lambda_-$. This proves
Assertion (i).

Similarly, Assertion (ii) of Proposition \ref{prop:matrixRiccati}
follows from the second claims of Lemmae \ref{lem4} and
\ref{lem5}. This concludes the proof of Proposition
\ref{prop:matrixRiccati}.  \cqfd

\section{Proof of the main result}
\label{sect:proofmain}

Let $X$ be as in the introduction, let $C$ be a closed strictly convex
subset of $X$ with nonempty interior, whose boundary $S=\partial C$ is
${\rm C}^{1,1}$-smooth, and let $\vec{n}$ be the inward normal unit
vector field along $S$. Note that $\vec{n}$ is well-defined (since $S$
is a ${\rm C}^{1}$ strictly convex hypersurface). It is locally
Lipschitz, hence is differentiable at almost every point of $S$ (for
the (well defined) Lebesgue measure on $S$).

Let ${\rm II}_S:TS\oplus TS\ra\RR$ be the (almost-everywhere defined)
scalar second fundamental form of $S$ associated to the inward normal
unit vector field $\vec{n}$ along $S$, that is, with
$\langle\cdot,\cdot\rangle$ the first fundamental form,
\begin{equation}\label{eq:definsecondformfond}
{\rm  II}_S(V,W)=\langle\nabla_VW,\vec{n}\rangle=
-\langle\nabla_V\vec{n},W\rangle\;,
\end{equation}
where $V,W$ are tangent vectors to $S$ at the same point, extended to
Lipschitz vector fields on a neighbourhood of this point, which are
tangent to $S$ at every point of $S$, and are differentiable at every
twice differentiable point of $S$. The definition of ${\rm II}_S$
depends on the choice between $\vec{n}$ and $-\vec{n}$, and the
various references differ on that point;
we have chosen the inward pointing vector field in order for the
symmetric bilinear form ${\rm II}_S$ to be nonnegative, by convexity
of $C$.

As a motivation, here is a short proof that if $C$ is
$\epsilon$-strictly convex, then $\overline{\rm II}_S\leq
b\coth(b\epsilon)$ almost everywhere.  For every $x$ in $S$, let
$y=\exp_x(\epsilon \,\vec{n}(x))$. Note that the sphere
$S_X(y,\epsilon)$ of center $y$ and radius $\epsilon$ in $X$ is
contained in $C$, as $C$ is $\epsilon$-strictly convex. Locally over
the tangent space $T_xS=T_x(S_X(y,\epsilon))\subset T_xX$, the graph
of $S_X(y,\epsilon)$ is above the graph of $S$ (when $\vec{n}$ points
upwards).  Hence, for almost every $x$ in $S$, for every $v$ in
$T_xS$, we have ${\rm II}_S(v,v)\leq {\rm II}_{S_X(y,\epsilon)}(v,v)$.
As the sectional curvature of $X$ is at least $-b^2$, we have by
comparison $\overline{\rm II}_{S_X(y,\epsilon)}\leq b\coth(b\epsilon)$
(see for instance \cite[page 175]{Pet}). 

\medskip
We now explain the curvature equation that will allow
us to apply the results of Section \ref{sec:compaRiccati}.

For every $t\in\RR$ and $x$ in $S$, let $x_t=\exp_x(t \,\vec{n}(x))$,
and let $N$ be the vector field defined by $N(x_t)=\dot{x}_t$ on an
open neighbourhood $U$ of $S$ (containing $X-C$ and $x_s$ for $s\in
[0,t]$ if it contains $x_t$). Note that $N$ is locally Lipschitz on
$U$, is differentiable at each $x_t\in U$ such that $S$ is twice
differentiable at $x$, and is smooth along each geodesic $t\mapsto
x_t$ while it stays in $U$. Identify the tangent space $T_{x_t}X$ with
$T_{x}X$ by the parallel transport $|\!|_{x_t}^{x}\colon T_{x_t}X\ra
T_{x}X$ along the geodesic $s\mapsto x_s$ whenever it stays in $U$.
Let $A(t)$ (which, for every twice differentiable point $x$ in $S$, is
defined and smooth at every $t$ such that $x_t\in U$) be the symmetric
endomorphism of $T_{x}S$ defined by $v\mapsto -|\!|_{x_t}^{x}
\big(\nabla_{|\!|^{x_t}_{x}(v)}N\big)$. Let $R(t)$ (which is defined
for every $x\in S$ and $t\in\RR$, is locally Lipschitz in $(x,t)$, and
is smooth in $t$) be the symmetric endomorphism of $T_{x}S$ defined by
$v\mapsto |\!|_{x_t}^{x} R(|\!|_{x}^{x_t}(v), \dot{x}_t) \dot{x}_t$,
where $R(\cdot,\cdot)\,\cdot\,$ is the curvature tensor of $X$.

For every twice differentiable point $x$ in $S$, the map $t\mapsto
A(t)$, defined and smooth on a neighbourhood of $0$, satisfies on it
the following matrix Riccati equation (see for instance \cite[page
44]{Pet}
\begin{equation}\label{eq:ricattiequation}
-\dot{A}(t) + A(t)^2 + R(t)=0\;.
\end{equation}

For every twice differentiable point $x\in S$ and every unit tangent
vector $v\in T^1_xS$, note that $\langle R(t)v,v\rangle$, being the
sectional curvature of the plane generated by the orthonormal tangent
vectors $|\!|_{x}^{x_t}v$ and $\dot{x}_t$ at $x_t$, is at most $-a^2$,
and at least $-b^2$. Hence, the map $t\mapsto A(t)$ satisfies the two
inequalities \eqref{eq:inegmatrixriccatti}.

\medskip
\noindent{\bf Proof of the first assertion of Theorem
  \ref{theo:maintheo}. } If $C$ is $\epsilon$-strictly convex, then
there exists a nonempty convex subset $C'$ of $X$ such that
$C=\N_\epsilon C'$. For every $t\in[0,\epsilon[$, the map $x\mapsto
x_t$ is a ${\rm C}^{1,1}$-diffeomorphism from $S$ to
$S_t=\partial\N_{\epsilon-t} C'$. For every twice differentiable point
$x$ of $S$, the point $x_t$ is a twice differentiable point of $S_t$,
and $A(t)$ is well-defined on $[0,\epsilon[$ (see \cite[\S 3]{Wal} for
these two facts). By the definition of the endomorphism $A(t)$, and
since the parallel transport preserves the first fundamental form, we
have for every $v$ in $T^1_{x_t}S_t$,
\begin{equation}\label{eq:relatIIstAt}
{\rm II}_{S_t}(v,v)=-\langle \nabla_vN,v\rangle=
\langle A(t)|\!|^{x_t}_{x}v,|\!|^{x_t}_{x}v\rangle \;.
\end{equation}
Since $S_t$ is locally convex, its second fundamental form is
nonnegative at each twice differentiable point, hence the endomorphism
$A(t)$ is nonnegative for $t\in[0,\epsilon[$, for almost every $x\in
S$.  The first assertion of Theorem \ref{theo:maintheo} now follows
from the second assertion of Proposition \ref{prop:matrixRiccati}.

\medskip
\noindent{\bf Proof of the second assertion of Theorem
  \ref{theo:maintheo}. } First assume that $S$ is smooth.

Recall (see \cite[page 222]{BC}) that the $\Sigma$-Jacobi fields for a
$\operatorname{C}^1$-smooth submanifold $\Sigma$ are the variations of
its normal geodesics. More precisely, for every $x\in S$, a map
$J\colon\RR\,\ra TX$ is a {\it $S$-Jacobi field} along the normal
geodesic $\tau\colon t\mapsto x_{t}$ to $S$ at $x$ if there exists a
$\operatorname{C}^1$-smooth map $f:\RR^2\ra X$ such that for every
$s,t\in\RR$, we have $f(t,0)=\tau(t)$, the map $t\mapsto f(t,s)$ is a
geodesic line in $X$ which starts at time $t=0$ perpendicularly to
$S$, and $J(t)=\frac{\partial f} {\partial s} (t,0)$. Note that by
Schwarz' theorem, the vector field $J$ commutes with the vector field
$N=\frac{\partial f} {\partial t} (t,0)$ along $\tau$. Since the
Riemannian connection is torsion-free, by the definition of $A(t)$ and
denoting again by $J(t)$ the parallel transport of $J(t)$ from
$T_{x_t}X$ to $T_xX$, we have
\begin{equation}\label{eq:SJacobi}
\dot J(t)=-A(t)J(t)\;.
\end{equation}

Recall that $t_0\geq 0$ is a {\it focal time} for $x_0\in S$ if the
differential of $x\mapsto x_{t_0}$ is not injective at $x_0$. Note
that $t_0$ is a focal time for $x\in S$ if and only if there is a
nonzero $S$-Jacobi field $J$ along $t\mapsto x_t$ which vanishes at
$x_{t_0}$, and that for every $\epsilon'>0$, by the triviality of the
normal bundle to $S$ and the convexity of $C$, the map $(x,t)\mapsto
x_t$ from $S\times ]-\infty,\epsilon']$ to $X$ is a proper immersion
if there is no focal time in $[0,\epsilon']$, see \cite[\S 11.3]{BC}.

To prove the second assertion of Theorem \ref{theo:maintheo}, let us
assume that
$$
b\tanh (b\epsilon)\le\underline{\rm II}_S
\le\overline{\rm II}_S\leq a\coth(a\epsilon).
$$
By the definition of the endomorphism $A(0)$, we have ${\rm II}_{S}
(v,v) = -\langle \nabla_v\vec{n},v\rangle= \langle A(0)v,v\rangle$ for
every $v$ in $T^1_{x}S$. Hence by the first assertion of Proposition
\ref{prop:matrixRiccati}, the endomorphism $A(t)$ is defined and
positive for all $t\in[0,\epsilon[\,$.

We claim that no nonzero $S$-Jacobi field along a normal geodesic
$t\mapsto x_t$ to $S$ vanishes in $[0,\epsilon[\,$. Indeed, since
$A(t)$ is defined for every $t\in [0,\epsilon[$, a $S$-Jacobi field,
which satisfies the first order linear equation \eqref{eq:SJacobi},
vanishes at one point of $[0,\epsilon[$ if and only if it vanishes at
all points of $[0,\epsilon[$. Since the biggest eigenvalue of $A(0)$
is at most $a\coth(a\epsilon)$, the claim also follows from the
comparison theorem for $S$-Jacobi fields \cite[Thm.~3.4]{Esch} applied
to $S$ and the sphere of radius $\epsilon$ in the real hyperbolic
space of constant curvature $-a^2$ (and keeping in mind that the sign
convention for $A$ in \cite{Esch} is different from ours).

Hence $x\mapsto x_t$ is a proper immersion of $S$, whose image we
denote by $S_t$, for every $t\in[0,\epsilon[\,$. By the definition of
the endomorphism $A(t)$, Equation \eqref{eq:relatIIstAt} is still
valid, and hence $S_t$ is, for every $t\in[0,\epsilon[\,$, a strictly
convex immersed hypersurface whose smallest eigenvalue of the second
fundamental form has a positive lower bound. Let $t_*\geq 0$ be the
upper bound of all $t\in[0,\epsilon[$ such that the map $(x,u)\mapsto
x_u$ from $S\times [0,t]$ to $X$ is an embedding.

If $t_*$ was strictly less than $\epsilon$, then there would exist a
sequence $(t_n)_{n\in\NN}$ in $[0,t_*[$ converging to $t_*$ and two
sequences $(x_n)_{n\in\NN}$ and $(y_n)_{n\in\NN}$ in $S_{t_n}$ such
that $d(x_n,y_n)$ tends to $0$ and $N(x_n)+|\!|_{y_n}^{x_n}N(y_n)$
converges to $0$. Hence as $n\ra+\infty$, the tangent planes at $x_n$
and $y_n$ are closer and closer, and the two germs of the hypersuface
$S_{t_n}$ at $x_n$ and $y_n$, being contained between them, are more
and more flat. But this contradicts the fact that the smallest
eigenvalue of the second fundamental form of $S_{t_*}$ has a positive
lower bound.

If $t_*=\epsilon$, then for every $s\in[0,\epsilon[\,$, the subset
$C_s= C- \bigcup_{t\in [0,s[} S_t$ is closed, and convex since its
boundary is locally convex, by the Schmidt theorem (see for
instance \cite[Appendix]{Hei} as explained in \cite[page
323]{Ale}). Define $C=\bigcap_{s\in[0,\epsilon[} C_s$, which is a
closed convex subset. We have $C=\N_\epsilon C'$ by construction, so
that $C$ is $\epsilon$-strictly convex. This ends the proof of the
second assertion of Theorem \ref{theo:maintheo} when $S$ is smooth.

\bigskip The following approximation result will allow us to extend
the result from the smooth case to the general case.

\bprop \label{prop:gw}
Let $a,b,\alpha,\beta>0$. Let $M$ be a complete simply
connected smooth Riemannian manifold of dimension $m\geq 2$ with
pinched sectional curvature $-b^2\leq K\leq -a^2<0$. Let $C$ be a
closed strictly convex subset of $M$ with nonempty interior, whose
boundary $S=\partial C$ is ${\rm C} ^{1,1}$-smooth and satisfies the
inequalities $\alpha\leq \underline{\rm II}_S\leq\overline{\rm II}_S\leq
\beta$ almost everywhere. Then for every $\eta>0$,
$\alpha'\in\;]0,\alpha[$ and $\beta'\in \; ]\beta,+\infty[$, there
exists a closed convex subset $C'$ in $X$ containing $C$, with smooth
boundary $S'=\partial C'$ such that $C'\subset \N_\eta C$ and
$\alpha'\leq \underline{\rm II}_{S'}\leq\overline{\rm II}_{S'}\leq
\beta'$.
\eprop

\dem 
We denote by $\exp:TM\ra M$ the Riemannian exponential map, by
$\parallel_{x}^y\,:T_xM\ra T_yM$ the parallel transport along the
geodesic from $x$ to $y$ in $M$, by $\nabla f:M\ra TM$ the Riemannian
gradient of a $\operatorname{C}^1$ map $f:M\ra\RR$, by $\pi:TM\ra M$
the canonical projection, by $TTM=V\oplus H$ the orthogonal
decomposition into the vertical and horizontal subbundles of $TTM\ra
TM$ defined by the Riemannian metric of $M$ (with $V_{v}=T_{\pi(v)}M$
and $T\pi_{\mid H_v}: H_v\ra T_{\pi(v)}M$ a linear isomorphism for
every $v\in TM$), and by $\pi_V:TTM\ra V$ the bundle projection to the
vertical factor parallel onto the horizontal one. Recall that the
covariant derivative of a $\operatorname{C}^1$ vector field $Y:M\ra
TM$ is defined by $\nabla Y=\pi_V\circ TY$. Also recall that if
$f:M\ra\RR$ and $g:M'\ra M$ are $\operatorname{C}^1$ maps, then
\begin{equation}\label{eq:gradientcomposition}
\nabla (f\circ g)=(Tg)^*\circ (\nabla f)\circ g\;,
\end{equation}
where $h^*:TM\ra TM'$ is the adjoint bundle morphism of a bundle
morphism $h:TM'\ra TM$ defined by the Riemannian metrics. Here is a
proof by lack of reference. 
For every $x\in M'$ and $Z\in T_xM'$, we
have
\begin{align*}
\langle\big(\nabla (f\circ g)\big) (x),Z\rangle_x&
=d(f\circ g)_x(Z)=df_{g(x)}(Tg(Z))
=\langle\nabla f ( g (x)),Tg(Z)\rangle_{g(x)}\\ &=
\langle(Tg)^*\big(\nabla f ( g (x))\big),Z\rangle_{x}\;.
\end{align*}

\medskip The main tool we use to prove Proposition \ref{prop:gw} is
the {\it Riemannian convolution smoothing} process of Greene and
Wu. Introduced in \cite[page 646]{GW1}, it has already been used for
instance in \cite[Theo.~2]{GW2} to approximate continuous strictly
convex functions on Riemannian manifolds by smooth ones. A new
property of this process introduced in this paper is a good control of
its second order derivatives. The smoothing is defined as follows. Let
$\psi:\RR\ra[0,+\infty[$ be a smooth map with compact support
contained in $[-1,1]$, constant on a neighbourhood of $0$, such that
$\int_{v\in\RR^m} \psi(\|v\|)\;d\lambda(v)= 1$, where $d\lambda$ is
the Lebesgue measure on the standard Euclidean space $\RR^m$.  For
every $\kappa>0$, let $\psi_\kappa :t\mapsto
\frac{1}{\kappa^m}\psi(\frac{t}{\kappa})$, whose support is contained
in $[-\kappa,\kappa]$. For every continuous map $f:M\ra \RR$, define a
map $f_\kappa:M\ra \RR$ by
$$
f_\kappa:x\mapsto \int_{v\in T_xM} \psi_\kappa(\|v\|)\; 
f(\exp_x v) \;d\lambda_x(v)\;,
$$
where $d\lambda_x$ is the Lebesgue measure on the Euclidean space
$T_xM$. Note that $f_\kappa$ is nonnegative if $f$ is nonnegative.

Given $\eta,\alpha',\beta'$ as in the statement of Proposition
\ref{prop:gw}, let us prove that if $f$ is the distance function to
$C$, for some $\kappa,t>0$ small enough, then
$C'={f_\kappa}^{-1}([0,t])$ satisfies the conclusions of Proposition
\ref{prop:gw}.

\medskip Let $j_x=\frac{d(\exp_x)_*\lambda_x}{d\operatorname{vol}}:
M\ra [0,+\infty]$ be the jacobian of the map $\exp_x$ from $T_{x}M$ to
$M$, with respect to the Lebesgue measure $d\lambda_{x}$ on $T_{x}M$
and the Riemannian measure $d\operatorname{vol}$ on $M$.  Since $M$ is
complete, simply connected and negatively curved, for every $x$ in
$M$, the map $\exp_x$ is a smooth diffeomorphism, whose inverse and
whose jacobian are hence well-defined, and depend smoothly on $(x,y)$
where $y$ is the variable point in $M$. For every continuous map
$f:M\ra \RR$, by an easy change of variables, we have
$f_\kappa:x\mapsto \int_{y\in M} \psi_\kappa (\|\exp_x^{-1}(y)\|)
\;f(y)\;j_x(y)\; d\operatorname{vol}(y)$. Since $\psi$ is constant
near $0$ and by a standard argument of differentiation under the
integral sign, the map $f_\kappa$ is smooth.

If $f:M\ra \RR$ is a $1$-Lipschitz map, then since $\int_{v\in T_xM}
\psi_\kappa(\|v\|) \; d\lambda_x(v)=1$, since $d(x,\exp_x v)=
\|v\|$ and since $v\mapsto \psi_\kappa(\|v\|)$ vanishes outside
$\{v\in T_xM\;:\;\|v\|\leq\kappa\}$, we have, for every $x\in M$,
\begin{equation}\label{eq:unifproche}
|f_\kappa(x)-f(x)|\leq \kappa\;.
\end{equation}

For every $x_0\in M$ and $v\in T_{x_0}M$, let $g_v:M\ra M$ be the map
$x\mapsto\exp_x(\parallel_{x_0}^xv)$. Note that $g_0$ is the identity
map, and that $(v,x)\mapsto g_v(x)$ is a smooth map. Since the
parallel transport is an isometry, for every continuous map
$f:M\ra \RR$, we have as in \cite{GW1}
$$
f_\kappa:x\mapsto\int_{v\in T_{x_0}M}\psi_\kappa(\|v\|)\;f\circ
g_v(x)\;d\lambda_{x_0}(v)\;.
$$
By Equation \eqref{eq:gradientcomposition}, and by differentiation
under the integral sign, if $f$ is $\operatorname{C}^1$ on
$B(x_0,2\kappa)$ and $x\in B(x_0,\kappa)$, we have
$$
\nabla f_\kappa(x)=\int_{v\in T_{x_0}M}\psi_\kappa(\|v\|)\;
(Tg_v)^*\circ (\nabla f)\circ g_v (x)\;d\lambda_{x_0}(v)\;.
$$

Let now $f:M\ra \;[0,+\infty[$ be the distance map to $C$, which is
$1$-Lipschitz on $M$ and is a ${\rm C}^{1,1}$ Riemannian submersion
outside $C$. Recall that (adapting \cite[\S 2.4.1]{Pet}
to the $\operatorname{C}^{1,1}$ regularity), the Hessian $Hf=\nabla^2
f : TM\ra TM$ of $f$ is almost everywhere defined. For every $t>0$,
the level hypersurface $S_t=f^{-1}(t)$ is $\operatorname{C}^{1}$ with
inward normal unit vector field equal to $-\nabla f$ along
$S_t$, and is twice differentiable at almost every point. By Equation
\eqref{eq:definsecondformfond}, the second fundamental form of $S_t$
at a twice differentiable point $x$ satisfies, for every $Z\in
T_xS_t$,
$$
{\rm II}_{S_t}(Z,Z)=\langle H f(Z) \,,\, Z\rangle\;.
$$
Let $\alpha''\in\;]\alpha',\alpha[$ and $\beta''\in\;]\beta,\beta'[
\,$. Let $\eta'\in\;]0,\eta]$ be small enough so that $\alpha'' \leq
\underline{\rm II}_{S_s} \leq\overline{\rm II}_{S_s} \leq \beta''$
almost everywhere on $S_s$ for every $s\in[0,\eta']$. Let
$t\in[\frac{\eta'}{3},\frac{2\eta'}{3}]$. Then by Equation
\eqref{eq:unifproche}, for every $\kappa\in\;]0,\frac{\eta'}{12}]$, we
have
$$
C\subset \N_{3\kappa}C\subset {f_\kappa}^{-1}([0,t])
\subset \N_{\eta'}C\subset \N_{\eta}C\;.
$$

By the linearity of $\pi_V$ and again by differentiation almost
everywhere of Lipschitz maps under the integral sign, if
$d(x_0,C)>2\kappa$ and $Z\in T_{x_0}M$, we have
$$
Hf_\kappa(Z)= \nabla^2 f_\kappa(Z)= \int_{v\in
  T_{x_0}X}\psi_\kappa(\|v\|)\;\pi_V\circ T\big((Tg_v)^*\big)\circ
T(\nabla f)\circ (Tg_v) (Z) \;d\lambda_{x_0}(v)\;.
$$
The maps $T\big((Tg_v)^*\big)$ and $Tg_v$ are the identity maps of
respectively $TTM$ and $TM$ when $v=0$, and they depend continuously
of $v$ for the uniform convergence of maps, since $M$ has pinched
curvature. If $\kappa$ is small enough, since $\pi_V$ is
$1$-Lipschitz, for every $x_0\in M$ such that $\frac{\eta'}{3}\leq
d(x_0,C)\leq\frac{2\eta'}{3}$, we hence have
$$ 
\alpha'\leq \min_{Z\in T^1_{x_0}M}\langle Hf_\kappa(Z),Z\rangle\leq 
\max_{Z\in T^1_{x_0}M}\langle Hf_\kappa(Z),Z\rangle \leq \beta'\;.
$$
Using Sard's theorem, let $t\in[\frac{\eta'}{3}, \frac{2\eta'}{3}]$ be
such that ${f_\kappa}^{-1}(t)$ is smooth, and we have
$$ 
\alpha'\leq \underline{\rm II}_{{f_\kappa}^{-1}(t)} \leq
\overline{\rm II}_{{f_\kappa}^{-1}(t)} \leq \beta'\;.
$$
In particular, the smooth hypersurface ${f_\kappa}^{-1}(t)$ is
strictly convex, since $\alpha'>0$.

Using again the Schmidt theorem, it is then easy to check that
$C'={f_\kappa}^{-1}([0,t])$ satisfies the conclusion of Proposition
\ref{prop:gw}.  
\cqfd

\medskip 
Now, if $C$ satisfies the assumption of the second assertion
of Theorem \ref {theo:maintheo}, for every $n\in\NN$ bigger than some
$N\in\NN$, there exists, by Proposition \ref{prop:gw}, a closed convex
subset $C_n$ with smooth boundary containing $C$ such that $C_n\subset
\N_{\frac{1}{n}}C$ and $b\tanh (b(\epsilon-\frac{1}{n})) \leq
\underline{\rm II}_S \leq \overline{\rm II}_S\leq
a\coth(a(\epsilon-\frac{1}{n}))$.  By the already proven smooth case
of the second assertion of Theorem \ref {theo:maintheo}, $C_n$ is
hence $(\epsilon-\frac{1}{n})$-strictly convex. Let $C'_n$ be a closed
convex subset such that $C_n=\N_{\epsilon-\frac{1}{n}}(C'_n)$.

Then $C$ is the closed $\epsilon$-neighbourhood of the intersection of
the $C'_n$'s. Indeed, since $\N_{\epsilon-\frac{1}{n}}(C'_n)$ is
contained in $\N_{\epsilon}(C'_n)$ for every $n\geq N$, the set
$C=\bigcap_{n\geq N} C_n$ is contained in
$\N_{\epsilon}\big(\bigcap_{n\geq N}C'_n\big)$.  Conversely, let $x\in
\N_{\epsilon}\big(\bigcap_{n\geq N}C'_n\big)$.  Since
$\N_{\epsilon}(C'_n)=\N_{\frac{1}{n}}C_n\subset \N_{\frac{2}{n}}C$,
for every $n\geq N$, there exists $x_n\in C$ such that $d(x,x_n)\leq
\frac{2}{n}$. By a compactness argument and since $C$ is closed, we
hence have $x\in C$.

This proves that $C$ is $\epsilon$-strictly convex, and concludes the
proof of Theorem \ref{theo:maintheo}.

\bigskip
{\small\noindent \begin{tabular}{l} 
Department of Mathematics and Statistics, P.O. Box 35\\ 
40014 University of Jyv\"askyl\"a, FINLAND.\\
{\it e-mail: parkkone@maths.jyu.fi}
\end{tabular}
\medskip

\noindent \begin{tabular}{l}
DMA, UMR 8553 CNRS\\
Ecole Normale Sup\'erieure, 45 rue d'Ulm\\
75230 PARIS Cedex 05, FRANCE\\
{\it e-mail: Frederic.Paulin@ens.fr}
\end{tabular} and \begin{tabular}{l}
D\'epartement de math\'ematique, B\^at.~425\\
Universit\'e Paris-Sud 11\\
91405 ORSAY Cedex, FRANCE\\
{\it e-mail: frederic.paulin@math.u-psud.fr}
\end{tabular}
}

\end{document}

%% file: macros_angl.tex
\usepackage[T1]{fontenc}
\usepackage[dvips]{color}
\usepackage{amssymb,amsmath,epsfig}
\usepackage{mathrsfs}
\renewcommand\mathcal{\mathscr}
\usepackage{comment}

\newtheorem{defi}{Definition}
\newtheorem{prop}[defi]{Proposition}
\newtheorem{theo}[defi]{Theorem}
\newtheorem{conj}[defi]{Conjecture}
\newtheorem{lemm}[defi]{Lemma}
\newtheorem{coro}[defi]{Corollary}
\newtheorem{rema}[defi]{Remark}
\newtheorem{exem}[defi]{Example}
\newtheorem{exems}[defi]{Examples}

\newcommand{\bdefi}{\begin{defi}}
\newcommand{\edefi}{\end{defi}}
\newcommand{\bprop}{\begin{prop}}
\newcommand{\eprop}{\end{prop}}
\newcommand{\btheo}{\begin{theo}}
\newcommand{\etheo}{\end{theo}}
\newcommand{\blemm}{\begin{lemm}}
\newcommand{\brema}{\begin{rema}}
\newcommand{\erema}{\end{rema}}
\newcommand{\bexer}{\begin{exem}}
\newcommand{\eexer}{\end{exem}}
\newcommand{\bexems}{\begin{exems}}
\newcommand{\eexems}{\end{exems}}
\newcommand{\bconj}{\begin{conj}}
\newcommand{\econj}{\end{conj}}
\newcommand{\elemm}{\end{lemm}}
\newcommand{\bcoro}{\begin{coro}}
\newcommand{\ecoro}{\end{coro}}
\newcommand{\dem}{\noindent{\bf Proof. }}

\newcommand{\A}{{\cal A}}

\newcommand{\N}{{\cal N}}

\newcommand{\maths}[1]{{\mathbb #1}}  

\newcommand{\RR}{\maths{R}}
\newcommand{\NN}{\maths{N}}

\newcommand{\ra}{\rightarrow}

\newcommand{\cqfd}{\hfill$\Box$}

\newcounter{fig}

